\documentclass[10pt]{amsart}

\usepackage{amssymb}

\newtheorem{thm1}{Theorem}[section]

\newtheorem{rem1}[thm1]{Remark}

\newtheorem{cor1}[thm1]{Corollary}

\newtheorem{prop1}[thm1]{Proposition}
\newtheorem{ex1}[thm1]{Example}

\begin{document}

\title[]
{Arithmetical rank and cohomological dimension of generalized binomial edge ideals}
\author[]{Anargyros Katsabekis}
\address { Department of Mathematics, University of Ioannina, 45110 Ioannina, Greece} \email{katsampekis@uoi.gr}

\begin{abstract} Let $G$ be a connected and simple graph on the vertex set $[n]$. To the graph $G$ one can associate the generalized binomial edge ideal $J_{m}(G)$ in the polynomial ring $R=K[x_{ij}: i \in [m], j \in [n]]$. We provide a lower bound for the cohomological dimension of $J_{m}(G)$. We also study when $J_{m}(G)$ is a cohomologically complete intersection. Finally, we show that the arithmetical rank of $J_{2}(G)$ equals the projective dimension of $R/J_{2}(G)$ in several cases.
\end{abstract}
\maketitle

\section{Introduction}

Let $R$ be a polynomial ring in finitely many variables over a field $K$ and let $I$ be an ideal of $R$. A classical problem in Algebraic Geometry is to determine the minimal number of equations needed to define set-theoretically the algebraic set corresponding to $I$. When $K$ is algebraically closed, it is equivalent to the problem of finding the smallest integer $s$ for which there exist polynomials $f_{1},\ldots,f_{s}$ in $I$ such that ${\rm rad}(I)={\rm rad}(f_{1},\ldots,f_{s})$. This integer is called the {\em arithmetical rank} of $I$ and will be denoted by ${\rm ara}(I)$. A lower bound for ${\rm ara}(I)$ is provided by the so-called {\em cohomological dimension} of $I$, which is defined as follows: ${\rm cd}(I)={\rm max}\{i \in \mathbb{Z}: H_{I}^{i}(R) \neq 0\}$, where $H_{I}^{i}(R)$ denotes the $i$-th local cohomology module of $R$ with support in $I$. In this paper we study the latter problem for a special class of polynomial ideals, namely the generalized binomial edge ideals of graphs.

Let $m \geq 2$, $n \geq 2$ be integers and $R=K[x_{ij}: i \in [m], j \in [n]]$  be the polynomial ring over a field $K$ in $mn$ variables $x_{ij}$, where $i=1,\ldots,m$ and $j=1,\ldots,n$. Let $G$ be a connected and simple graph on the vertex set $[n]$. The {\em generalized binomial edge idea}l $J_{m}(G)$ of $G$ is generated by all the $2$-minors of the $m \times n$ matrix $X=(x_{ij})$ of the form $[k,l|i,j]=x_{ki}x_{lj}-x_{kj}x_{li}$, where $1 \leq k<l \leq m$ and $\{i,j\}$ is an edge of $G$ with $i<j$. When $m=2$, $J_{m}(G)$ coincides with the classical binomial edge ideal introduced independently in \cite{HHHK} and \cite{Oht}.

 By \cite[Theorem 2]{Rauh}, $J_{m}(G)$ has a squarefree initial ideal with respect to a term order $<$. Thus $J_{m}(G)$ is a radical ideal. By \cite[Corollary 2.7]{CV}, ${\rm depth}(R/J_{m}(G))={\rm depth}(R/{\rm in}_{<}(J_{m}(G)))$, and therefore ${\rm pd}(R/J_{m}(G))={\rm pd}(R/{\rm in}_{<}(J_{m}(G)))$. Since ${\rm in}_{<}(J_{m}(G))$ is generated by squarefree monomials, its cohomological dimension equals the projective dimension (see \cite{Lyu}). So ${\rm pd}(R/J_{m}(G))={\rm cd}({\rm in}_{<}(J_{m}(G)))$. By \cite[Proposition 3.6]{CV}, the inequality ${\rm pd}(R/J_{m}(G)) \leq {\rm cd}(J_{m}(G))$ holds.

For a generalized binomial edge ideal $J_{m}(G)$ the following inequalities hold: $${\rm ht}(J_{m}(G)) \leq {\rm pd}(R/J_{m}(G)) \leq {\rm cd}(J_{m}(G)) \leq {\rm ara}(J_{m}(G)) \leq \mu(J_{m}(G)).$$ Here ${\rm ht}(J_{m}(G))$ denotes the height and $\mu(J_{m}(G))$ denotes the minimal number of generators of $J_{m}(G)$. The ideal $J_{m}(G)$ is a {\em set-theoretic complete intersection} if ${\rm ara}(J_{m}(G))={\rm ht}(J_{m}(G))$, while $J_{m}(G)$ is a {\em cohomologically complete intersection} if ${\rm cd}(J_{m}(G))={\rm ht}(J_{m}(G))$. Clearly if $J_{m}(G)$ is a set-theoretic complete intersection, then it is also a cohomologically complete intersection. It is worth mentioning that ${\rm ara}(J_{m}(G))=mn-3$ when $G$ is the complete graph on $[n]$, see \cite[Theorem 2]{BS}. It is still an open problem to compute ${\rm ara}(J_{m}(G))$ when $G$ is a non-complete graph.

In Section 2 we determine all graphs $G$ for which either $\mu(J_{m}(G))={\rm ht}(J_{m}(G))$ or $\mu(J_{m}(G))={\rm ht}(J_{m}(G))+1$, see Theorems \ref{complete}, \ref{almost}. Additionally we provide (see Theorem \ref{intersection}) a lower bound for the cohomological dimension of the ideal $P_{\emptyset}(G) \cap P_{T}(G)$, where $P_{\emptyset}(G)$ and $P_{T}(G)$ are distinct minimal primes of $J_{m}(G)$. We also show that, for a non-complete graph $G$ with vertex connectivity $k(G)$, ${\rm cd}(J_{m}(G)) \geq mn-m-n+k(G)$.

In Section 3 we prove that $J_{m}(G)$ is a cohomologically complete intersection for $m \geq 3$ if and only if $G$ is the complete graph on $[n]$ and $K$ has positive characteristic. As a consequence, we derive that $J_{m}(G)$ is not a set-theoretic complete intersection for every $m \geq 3$. In fact there can be a large deviation between ${\rm ara}(J_{m}(G))$ and ${\rm ht}(J_{m}(G))$, see Proposition \ref{deviation}. We also study the case $m=2$, namely $J_{2}(G)$ is a binomial edge ideal. When $\mu(J_{2}(G))={\rm ht}(J_{2}(G))+1$ and $G$ is not a triangle, we show that $J_{2}(G)$ is a cohomologically complete intersection if and only if $J_{2}(G)$ is Cohen-Macaulay. Finally we provide classes of graphs such that ${\rm ara}(J_{2}(G))={\rm pd}(R/J_{2}(G))$, see Theorems \ref{11}, \ref{12}, Corollary \ref{2}, Remark \ref{almostbasic}, and Theorems \ref{even}, \ref{An}. 

The computations of this paper are performed by using CoCoA \cite{cocoa}.\\

\section{Cohomological dimension}

Let $G$ be a connected graph on the vertex set $V(G)=[n]$ and let $T \subset [n]$. By $G \setminus T$ we denote the graph that results from deleting all vertices in $T$ and their incident edges from $G$. Let $c(T)$ be the number of connected components of $G \setminus T$. We say that $T$ is a {\em cut set} of $G$ if either $T=\emptyset$ or $c(T \setminus \{i\})<c(T)$ for every $i \in T$. We denote by $\mathcal{C}(G)$ the collection of cut sets of $G$.

Let $G_{1},\ldots,G_{c(T)}$ denote the connected components of $G \setminus T$ and also let $\tilde{G}_{i}$ denote the complete graph on the vertices of $G_i$. Consider the ideal $$P_{T}(G)=\left(\bigcup_{i \in T}\{x_{1i},\ldots,x_{mi}\} \right)+J_{m}(\tilde{G}_{1})+\cdots+J_{m}(\tilde{G}_{c(T)})$$ of the polynomial ring $R$. Then $P_{T}(G)$ is a prime ideal for every $T \subset [n]$. Furthermore $P_{T}(G)$ is a minimal prime of $J_{m}(G)$ if and only if $T \in \mathcal{C}(G)$. By \cite[Theorem 7]{Rauh}, $J_{m}(G)=\cap_{T \in \mathcal{C}(G)} P_{T}(G)$. Since ${\rm ht}(P_{\emptyset}(G))=(m-1)(n-1)$, it holds that ${\rm dim}(R/J_{m}(G)) \geq m+n-1$.

\begin{rem1} \label{BaseRem} {\rm (1) If $J_{m}(G)$ is Cohen-Macaulay, then ${\rm dim}(R/J_{m}(G)) =m+n-1$.\\
(2) If for a set $T \subsetneqq [n]$ we have that $c(T)=1$, then $P_{\emptyset}(G) \subset P_{T}(G)$ and therefore $P_{T}(G)$ is not a minimal prime of $J_{m}(G)$. So, if $P_{T}(G)$ is a minimal prime of $J_{m}(G)$, then either $T=\emptyset$ or $c(T) \geq 2$.}
\end{rem1}

\begin{prop1} Let $G$ be a connected graph on the vertex set $[n]$. Then $${\rm ht}(J_{m}(G))={\rm min}\{(m-1)(n-c(T))+|T|: T \in \mathcal{C}(G)\}.$$
\end{prop1}

\noindent \textbf{Proof.}  Let $T \in \mathcal{C}(G)$ and $|V(\tilde{G}_{i})|=n_{i}$ for $1 \leq i \leq c(T)$. We have that $${\rm ht}(P_{T}(G))={\rm ht} \left(\bigcup_{i \in T}\{x_{1i},\ldots,x_{mi}\} \right)+\sum_{i=1}^{c(T)}{\rm ht}(J_{m}(\tilde{G}_{i}))=$$ $$m|T|+\sum_{i=1}^{c(T)} (m-1)(n_{i}-1)=m|T|+(m-1) \sum_{i=1}^{c(T)} (n_{i}-1).$$  It holds that $$\sum_{i=1}^{c(T)} (n_{i}-1)=\sum_{i=1}^{c(T)}n_{i}-\sum_{i=1}^{c(T)}1=n-|T|-c(T),$$ so ${\rm ht}(P_{T}(G))=(m-1)(n-c(T))+|T|$ and therefore $${\rm ht}(J_{m}(G))={\rm min}\{(m-1)(n-c(T))+|T|: T \in \mathcal{C}(G)\}. \ \ \square$$

We say that $J_{m}(G)$ is a {\em complete intersection} if $\mu(J_{m}(G))={\rm ht}(J_{m}(G))$, while $J_{m}(G)$ is an {\em almost complete intersection} if $\mu(J_{m}(G))={\rm ht}(J_{m}(G))+1$.

\begin{thm1} \label{complete} Let $G$ be a connected graph on the vertex set $[n]$. Then $J_{m}(G)$ is a complete intersection if and only if $m=2$ and $G$ is a path graph.
\end{thm1}

\noindent \textbf{Proof.} Since $G$ is a connected graph, it has at least $n-1$ edges and therefore $\mu(J_{m}(G)) \geq \frac{m(m-1)}{2} (n-1)$. If $J_{m}(G)$ is a complete intersection, then it is unmixed. Since ${\rm ht}(P_{\emptyset}(G))=(m-1) (n-1)$, we get that ${\rm ht}(J_{m}(G))=(m-1) (n-1)$. Thus $\mu(J_{m}(G))=(m-1)(n-1)$, so $(m-1) \geq \frac{m(m-1)}{2}$ and therefore $m \leq 2$. Consequently $m=2$. By \cite[Theorem 1]{Rin}, $J_{m}(G)$ is a complete intersection if and only if $G$ is a path graph. \hfill $\square$

\begin{thm1} \label{almost} Let $G$ be a connected graph on the vertex set $[n]$, where $n \geq 3$. Then $J_{m}(G)$ is an almost complete intersection if and only if $m=2$ and either \begin{enumerate} \item $G$ is a tree which is not a path obtained by adding an edge between two vertices of two paths or
\item $G$ is not a tree obtained by adding an edge between two vertices of a path or by attaching a path to each vertex of a triangle. 
\end{enumerate}
\end{thm1}

\noindent \textbf{Proof.} The graph $G$ is connected, so it has at least $n-1$ edges, thus $\mu(J_{m}(G)) \geq \frac{m(m-1)}{2} (n-1)$. If $J_{m}(G)$ is an almost complete intersection, then $\mu(J_{m}(G))={\rm ht}(J_{m}(G))+1$ and therefore ${\rm ht}(J_{m}(G)) \geq \frac{m(m-1)}{2} (n-1)-1$. But ${\rm ht}(J_{m}(G)) \leq (m-1)(n-1)$, so $\frac{m(m-1)}{2} (n-1)-1 \leq (m-1)(n-1)$. Thus $(m-1)(n-1)(m-2) \leq 2$ from which we deduce that $m \leq 2$, and therefore $m=2$. By \cite[Theorems 4.3 and 4.4]{JKS}, $J_{m}(G)$ is an almost complete intersection if and only if $m=2$ and either \begin{enumerate} \item $G$ is a tree which is not a path obtained by adding an edge between two vertices of two paths or
\item $G$ is not a tree obtained by adding an edge between two vertices of a path or by attaching a path to each vertex of a triangle. \hfill $\square$
\end{enumerate}

\begin{rem1} {\em Let $G$ be a connected graph on the vertex set $[2]$, namely $G$ consists of the edge $\{1,2\}$. We have that $\mu(J_{m}(G))=\frac{m(m-1)}{2}$ and also ${\rm ht}(J_{m}(G))=m-1$. For $m=3$ the ideal $J_{m}(G)$ is an almost complete intersection, while $\mu(J_{m}(G))>{\rm ht}(J_{m}(G))+1$ for every integer $m>3$.}
\end{rem1}

Given a connected graph $G$ on the vertex set $[n]$, we say that $G$ is $k$-{\em vertex-connected} if $k<n$ and $G \setminus T$ is connected for every subset $T$ of $[n]$ with $|T|<k$. The {\em vertex connectivity} of $G$ is the maximum integer $k$ such that $G$ is $k$-vertex-connected. We shall denote by $k(G)$ the vertex connectivity of the graph $G$.

\begin{prop1} \label{Toobasic} Let $G$ be a connected graph on the vertex set $[n]$ and $P_{T}(G)$ be a minimal prime of $J_{m}(G)$, where $T \neq \emptyset$. Then ${\rm ht}(P_{\emptyset}(G)+P_{T}(G)) \geq (m-1)(n-1)+k(G)$.
\end{prop1}

\noindent \textbf{Proof.} The sum $P_{\emptyset}(G)+P_{T}(G)$ is generated by $\{x_{ki}x_{lj}-x_{kj}x_{li}|1 \leq k<l \leq m \ \textrm{and} \ i,j \in [n] \setminus T\} \cup \{x_{1i},\ldots,x_{mi}| i \in T\}$, so $${\rm ht}(P_{\emptyset}(G)+P_{T}(G))=(m-1)(n-|T|-1)+m|T|=(m-1)(n-1)+|T|.$$ It remains to show that $|T| \geq k(G)$. Note that $|T| \geq 1$. So if $k(G)=1$, then clearly $|T| \geq k(G)$. Suppose that $k(G) \geq 2$. By Remark \ref{BaseRem}, $c(T) \geq 2$ and therefore the graph $G \setminus T$ is disconnected. If $|T|<k(G)$, then the graph $G \setminus T$ is connected, a contradiction. Thus $|T| \geq k(G)$. \hfill $\square$

\begin{thm1} \label{intersection} Let $K$ be a field of any characteristic and $P_{T}(G)$ be a minimal prime of $J_{m}(G)$, where $T \neq \emptyset$. Then $${\rm cd}(P_{\emptyset}(G) \cap P_{T}(G)) \geq mn-m-n+|T|.$$
\end{thm1}

\noindent \textbf{Proof.} Recall that ${\rm ht}(P_{\emptyset}(G)+P_{T}(G))=(m-1)(n-1)+|T|$. Since $|T| \geq 1$, we get ${\rm ht}(P_{\emptyset}(G))<{\rm ht}(P_{\emptyset}(G)+P_{T}(G))$. Also ${\rm ht}(P_{T}(G))=(m-1)(n-c(T))+|T|$ and $c(T) \geq 2$, so ${\rm ht}(P_{T}(G))<{\rm ht}(P_{\emptyset}(G)+P_{T}(G))$. By \cite[Proposition 1.1.2]{Varbaro}, ${\rm cd}(P_{\emptyset}(G) \cap P_{T}(G)) \geq {\rm ht}(P_{\emptyset}(G)+P_{T}(G))-1$ and therefore ${\rm cd}(P_{\emptyset}(G) \cap P_{T}(G)) \geq (m-1)(n-1)+|T|-1$. Thus ${\rm cd}(P_{\emptyset}(G) \cap P_{T}(G)) \geq mn-m-n+|T|$. \hfill $\square$

\begin{ex1} \label{mimimal} {\rm Let $m=2$ and $G$ be a cycle with the edge set $$E(G)=\{\{1,2\},\{2,3\},\{3,4\},\{1,4\}\}.$$ Let $T=\{2,4\}$, then the ideal $P_{T}(G)=(x_{12},x_{22},x_{14},x_{24})$ is a minimal prime of $J_{2}(G)$. Also $P_{\emptyset}(G) \cap P_{T}(G)=J_{2}(G)+(x_{12}x_{24}-x_{14}x_{22})$. By Theorem \ref{intersection},   ${\rm cd}(P_{\emptyset}(G) \cap P_{T}(G)) \geq 4$. Let $$I=(x_{11}x_{22}-x_{12}x_{21}+x_{13}x_{24}-x_{14}x_{23},x_{12}x_{23}-x_{13}x_{22},x_{11}x_{24}-x_{14}x_{21},x_{12}x_{24}-x_{14}x_{22}).$$ Then $P_{\emptyset}(G) \cap P_{T}(G)={\rm rad}(I)$ since $(x_{11}x_{22}-x_{12}x_{21})^2 \in I$ and $(x_{13}x_{24}-x_{14}x_{23})^2 \in I$. So ${\rm ara}(P_{\emptyset}(G) \cap P_{T}(G)) \leq 4$ and therefore ${\rm cd}(P_{\emptyset}(G) \cap P_{T}(G)) \leq 4$. Thus ${\rm cd}(P_{\emptyset}(G) \cap P_{T}(G))=4$ and also ${\rm ara}(P_{\emptyset}(G) \cap P_{T}(G))=4$.}
\end{ex1}

Let $\mathcal{K}_{n}$ denotes the complete graph on the vertex set $[n]$.

\begin{rem1} \label{CDcomplete} {\rm It holds that \[ {\rm cd}(J_{m}(\mathcal{K}_{n}))=\begin{cases} 
      (m-1)(n-1), & {\rm if} \ \ {\rm char}(K)>0 \\
      mn-3, & {\rm if} \ \ {\rm char}(K)=0.
   \end{cases}
\] }
\end{rem1}

\begin{prop1}  Let $G$ be a path graph on the vertex set $[n]$, where $n \geq 3$, with the edge set $E(G)=\{\{1,2\},\{2,3\},\ldots,\{n-1,n\}\}$ and $T=\{2,3,\ldots,n-1\}$. \begin{enumerate} \item If $K$ is a field of positive characteristic, then ${\rm cd}(P_{\emptyset}(G) \cap P_{T}(G))=mn-m-2$. \item If $K$ is a field of characteristic $0$, then $mn-m-2 \leq {\rm cd}(P_{\emptyset}(G) \cap P_{T}(G)) \leq mn-3$. 
\end{enumerate}
\end{prop1}
\noindent \textbf{Proof.} Let $I=P_{\emptyset}(G) \cap P_{T}(G)$. Since $|T|=n-2$, we get ${\rm cd}(I) \geq m(n-1)-2$. Notice that $P_{T}(G)=\left(\cup_{i=2}^{n-1} \{x_{1i},\ldots,x_{mi}\} \right)$. Let $M=P_{\emptyset}(G)+P_{T}(G)$, then $M=P_{T}(G)+J_{m}(\mathcal{K}_{2})$, where $\mathcal{K}_{2}$ is the complete graph on the vertex set $\{1,n\}$. By \cite[Lemma 2.4]{BCM}, ${\rm cd}(M)=m(n-2)+{\rm cd}(J_{m}(\mathcal{K}_{2}))$. We distinguish the following cases: (1) ${\rm char}(K)>0$. Since ${\rm cd}(P_{\emptyset}(G))=(m-1)(n-1)$, ${\rm cd}(M)=m(n-1)-1$ and ${\rm cd}(P_{T}(G))=m(n-2)$, we deduce from the Mayer-Vietoris sequence
$$\cdots \rightarrow H_{M}^{i}(R) \rightarrow H_{P_{\emptyset}(G)}^{i}(R) \oplus H_{P_{T}(G)}^{i}(R) \rightarrow H_{I}^{i}(R) \rightarrow H_{M}^{i+1}(R) \rightarrow \cdots$$ that ${\rm cd}(I) \leq m(n-1)-2$ and therefore ${\rm cd}(I)=m(n-1)-2$.\\
(2) ${\rm char}(K)=0$. Since ${\rm cd}(P_{\emptyset}(G))=mn-3$, ${\rm cd}(M)=mn-3$ and ${\rm cd}(P_{T}(G))=m(n-2)$, we deduce from the Mayer-Vietoris sequence
$$\cdots \rightarrow H_{M}^{i}(R) \rightarrow H_{P_{\emptyset}(G)}^{i}(R) \oplus H_{P_{T}(G)}^{i}(R) \rightarrow H_{I}^{i}(R) \rightarrow H_{M}^{i+1}(R) \rightarrow \cdots$$ that $ {\rm cd}(I) \leq mn-3$. \hfill $\square$

\begin{thm1} \label{Basic1} Let $G$ be a connected graph on the vertex set $[n]$ and $K$ be a field of any characteristic. If $G$ is not the complete graph, then $${\rm cd}(J_{m}(G)) \geq mn-m-n+k(G).$$
\end{thm1}

\noindent \textbf{Proof.} Let $\mathcal{C}(G)=\emptyset \cup \{T_{1},\ldots,T_{s}\}$ be the cut sets of $G$, where $T_{i} \neq \emptyset$ for every $1 \leq i \leq s$. Since $G$ is not the complete graph, the ideal $J_{m}(G)$ is not prime. We can write $J_{m}(G)=P_{\emptyset}(G) \cap M$, where $M=\cap_{i=1}^{s} P_{T_i}(G)$ . It holds that ${\rm ht}(M) \leq {\rm ht}(P_{T_i}(G))$ and also ${\rm ht}(P_{T_i}(G))<{\rm ht}(P_{\emptyset}(G)+P_{T_{i}}(G))$, for every $1 \leq i \leq s$. There is $j \in \{1,\ldots,s\}$ such that ${\rm ht}(P_{\emptyset}(G)+P_{T_{j}}(G)) \leq {\rm ht}(P_{\emptyset}(G)+M)$. Consequently ${\rm ht}(M)<{\rm ht}(P_{\emptyset}(G)+M)$. Also ${\rm ht}(P_{\emptyset}(G)+P_{T_{j}}(G))>{\rm ht}(P_{\emptyset}(G))$, so ${\rm ht}(P_{\emptyset}(G)+M)>{\rm ht}(P_{\emptyset}(G))$. By \cite[Proposition 1.1.2]{Varbaro}, ${\rm cd}(J_{m}(G)) \geq {\rm ht}(P_{\emptyset}(G)+M)-1$. From Proposition \ref{Toobasic} we have that ${\rm ht}(P_{\emptyset}(G)+P_{T_{j}}(G)) \geq (m-1)(n-1)+k(G)$, so ${\rm cd}(J_{m}(G)) \geq mn-m-n+k(G)$. \hfill $\square$\\

A complete bipartite graph of the form $\mathcal{K}_{1,n-1}$ is a {\em star graph} with $n$ vertices.

\begin{prop1} \label{Pathcor} Let $G$ be a star graph on the vertex set $[n]$, where $n \geq 3$. \begin{enumerate} \item If $K$ is a field of positive characteristic, then ${\rm cd}(J_{m}(G))=(m-1)(n-1)$. \item If $K$ is a field of characteristic $0$, then $$(m-1)(n-1) \leq {\rm cd}(J_{m}(G)) \leq mn-3.$$
\end{enumerate}
\end{prop1}

\noindent \textbf{Proof.} The vertex connectivity of $G$ is $k(G)=1$, so ${\rm cd}(J_{m}(G)) \geq (m-1)(n-1)$. Let $i$ be the vertex of $G$ with degree $n-1$ and $T=\{i\}$. Then $J_{m}(G)=P_{\emptyset}(G) \cap P_{T}(G)$, where $P_{T}(G)=(x_{1i},x_{2i},\ldots,x_{mi})$. Let $M=P_{\emptyset}(G)+P_{T}(G)$, then $M=(x_{1i},x_{2i},\ldots,x_{mi})+J_{m}(\mathcal{K}_{n-1})$, where $\mathcal{K}_{n-1}$ is the complete graph on the vertex set $[n] \setminus \{i\}$. By \cite[Lemma 2.4]{BCM}, ${\rm cd}(M)=m+{\rm cd}(J_{m}(\mathcal{K}_{n-1}))$. So \[ {\rm cd}(M)=\begin{cases} 
      (m-1)(n-1)+1, & {\rm if} \ \ {\rm char}(K)>0 \\
      mn-3, & {\rm if} \ \ {\rm char}(K)=0.
   \end{cases}
\]
We distinguish the following cases: (1) ${\rm char}(K)>0$. Since ${\rm cd}(P_{\emptyset}(G))=(m-1)(n-1)$, ${\rm cd}(M)=(m-1)(n-1)+1$ and ${\rm cd}(P_{T}(G))=m$, we deduce from the Mayer-Vietoris sequence
$$\cdots \rightarrow H_{M}^{i}(R) \rightarrow H_{P_{\emptyset}(G)}^{i}(R) \oplus H_{P_{T}(G)}^{i}(R) \rightarrow H_{J_{m}(G)}^{i}(R) \rightarrow H_{M}^{i+1}(R) \rightarrow \cdots$$ that $ {\rm cd}(J_{m}(G)) \leq (m-1)(n-1)$, and therefore $ {\rm cd}(J_{m}(G))= (m-1)(n-1)$.\\
(2) ${\rm char}(K)=0$. Since ${\rm cd}(P_{\emptyset}(G))=mn-3$, ${\rm cd}(M)=mn-3$ and ${\rm cd}(P_{T}(G))=m$, we deduce from the Mayer-Vietoris sequence
$$\cdots \rightarrow H_{M}^{i}(R) \rightarrow H_{P_{\emptyset}(G)}^{i}(R) \oplus H_{P_{T}(G)}^{i}(R) \rightarrow H_{J_{m}(G)}^{i}(R) \rightarrow H_{M}^{i+1}(R) \rightarrow \cdots$$ that $ {\rm cd}(J_{m}(G)) \leq mn-3$. \hfill $\square$\\

\section{Arithmetical rank}

In this section first we determine all graphs $G$ for which $J_{m}(G)$ is a cohomologically complete intersection for $m \geq 3$.

\begin{prop1} \label{BasicProp} Let $G$ be a  connected graph on the vertex set $[n]$. If $J_{m}(G)$ is a cohomologically complete intersection, then $J_{m}(G)$ is Cohen-Macaulay.
\end{prop1}

\noindent \textbf{Proof.} Suppose that ${\rm cd}(J_{m}(G))={\rm ht}(J_{m}(G))$, so ${\rm pd}(R/J_{m}(G))={\rm ht}(J_{m}(G))$. Using the Auslander-Buchsbaum formula we get ${\rm depth}(R/J_{m}(G))={\rm dim} (R/J_{m}(G))$, thus $J_{m}(G)$ is Cohen-Macaulay. \hfill $\square$

\begin{rem1} \label{BasicRemark} {\rm If $J_{m}(G)$ is a set-theoretic complete intersection, then $J_{m}(G)$ is also Cohen-Macaulay.}
\end{rem1}

\begin{thm1} \label{cohbasic} Let $G$ be a connected graph on the vertex set $[n]$ and $m \geq 3$ be an integer. Then $J_{m}(G)$ is a cohomologically complete intersection if and only if $G$ is the complete graph and ${\rm char}(K)>0$.
\end{thm1}

\noindent \textbf{Proof.} If $G$ is the complete graph on the vertex set $[n]$ and ${\rm char}(K)>0$, then ${\rm cd}(J_{m}(G))={\rm ht}(J_{m}(G))=(m-1)(n-1)$. Conversely, if $J_{m}(G)$ is a cohomologically complete intersection, then $J_{m}(G)$ is Cohen-Macaulay from Proposition \ref{BasicProp}. By \cite[Corollary 4.3]{ACR}, $G$ is the complete graph on $[n]$. \hfill $\square$

\begin{thm1} Let $G$ be a connected graph on the vertex set $[n]$ and $m \geq 3$ be an integer. Then the ideal $J_{m}(G)$ is not a set-theoretic complete intersection.
\end{thm1}
\noindent \textbf{Proof.} Suppose that $J_{m}(G)$ is a set-theoretic complete intersection. By Theorem \ref{cohbasic}, $G$ is the complete graph on the vertex set $[n]$ and ${\rm char}(K)>0$. Notice that ${\rm ht}(J_{m}(G))=(m-1)(n-1)$. Thus ${\rm ara}(J_{m}(G))=mn-3>{\rm ht}(J_{m}(G))$ a contradiction to the fact that $J_{m}(G)$ is a set-theoretic complete intersection. \hfill $\square$\\

A {\em maximum clique} of a graph $G$ is a complete subgraph of $G$ that includes the largest possible number of vertices.

\begin{prop1}\label{Basic2} Let $G$ be a connected graph on the vertex set $[n]$ with $q$ edges and let $\mathcal{K}_{r}$ be a maximum clique of $G$, where $r \geq 2$. Then $${\rm ara}(J_{m}(G)) \leq \frac{m(m-1)}{2} \left(q-\frac{r(r-1)}{2} \right)+rm-3.$$
\end{prop1}

\noindent \textbf{Proof.} We have that $J_{m}(G)=J_{m}(\mathcal{K}_{r})+J_{m}(G')$, where $G'$ is a graph with $q-\frac{r(r-1)}{2}$ edges. Since  ${\rm ara}(J_{m}(\mathcal{K}_{r}))=rm-3$, we have that $${\rm ara}(J_{m}(G)) \leq \frac{m(m-1)}{2} \left(q-\frac{r(r-1)}{2} \right)+rm-3. \ \  \square$$

\begin{prop1} \label{deviation} Let $G$ be a star graph on the vertex set $[3]$ and $r \geq 1$ be an integer. Then for every integer $m \geq r+3$ it holds that $${\rm ara}(J_{m}(G))>{\rm ht}(J_{m}(G))+r.$$
\end{prop1}

\noindent \textbf{Proof.} By Proposition \ref{Pathcor}, ${\rm cd}(J_{m}(G)) \geq 2m-2$ and therefore ${\rm ara}(J_{m}(G)) \geq 2m-2$. Also ${\rm ht}(J_{m}(G))=m$. For every $m \geq r+3$ we have that $2m-2>m+r$, so ${\rm ara}(J_{m}(G))>{\rm ht}(J_{m}(G))+r$. \hfill $\square$\\

Next we prove that in various cases ${\rm ara}(J_{2}(G))={\rm pd}(R/J_{2}(G))$.\\ Let $G_1=(V(G_{1}),E(G_{1}))$, $G_2=(V(G_{2}),E(G_{2}))$ be graphs such that $G_{1} \cap G_{2}=\mathcal{K}_{m}$ is a complete graph, where $G_{1} \neq \mathcal{K}_{m}$ and $G_{2} \neq \mathcal{K}_{m}$. The new graph $G$ with the vertex set $V(G)=V(G_1) \cup V(G_2)$ and edge set $E(G)=E(G_{1}) \cup E(G_{2})$ is called the {\em clique sum} of $G_1$ and $G_2$ along $\mathcal{K}_{m}$, denoted by $G_{1} \cup_{\mathcal{K}_{m}} G_{2}$, see also \cite{JKS1}. If $m=1$, namely $\mathcal{K}_{1}=\{v\}$, then the clique sum of $G_1$ and $G_2$ along $v$ is denoted by $G_{1} \cup_{v} G_{2}$.

The {\em diamond graph} is the graph obtained from $\mathcal{K}_4$ by deleting one edge. Let $$f_{i,j}:=x_{1i}x_{2j}-x_{1j}x_{2i}.$$

\begin{thm1}  \label{11} Let $G$ be a connected graph on the vertex set $[n]$ with $n \geq 6$. Suppose that $G$ has a diamond subgraph $D$ and also $$G=D \cup_{v_1} T_{1} \cup_{v_2} \cdots \cup_{v_s}T_{s},$$ where $\{v_{1},\ldots,v_{s}\} \subset V(D)$, $s \geq 2$, $v_{h}$ are pairwise distinct and $T_{h}$ are trees. Then ${\rm ara}(J_{2}(G))={\rm pd}(R/J_{2}(G))=n-1$.
\end{thm1}
\noindent \textbf{Proof.} By \cite[Theorems 3.19 and 3.20]{BB}, ${\rm pd}(R/J_{2}(G)) \geq n-1$ and therefore ${\rm ara}(J_{2}(G)) \geq n-1$. Let $r_{h}$ be the number of edges of the tree $T_{h}$ for $1 \leq h \leq s$. The graph $G$ has $n=4+\sum_{h=1}^{s}r_{h}$ vertices. Let $\{v_{1},v_{5}\}$ be an edge of $T_{1}$ and $\{v_{2},v_{6}\}$ be an edge of $T_{2}$. Consider the subgraph $G_{1}=D \cup \{v_{1},v_{5}\} \cup \{v_{2},v_{6}\}$ of $G$. We will show that ${\rm ara}(J_{2}(G_{1})) \leq 5$. We distinguish the following cases:\\ (1) $E(D)=\{\{v_{1},v_{2}\},\{v_{2},v_{3}\},\{v_{3},v_{4}\},\{v_{1},v_{4}\},\{v_{1},v_{3}\}\}$. Let $$I=(f_{v_{1},v_{3}}+f_{v_{2},v_{6}},f_{v_{1},v_{5}}+f_{v_{3},v_{4}},f_{v_{1},v_{2}},f_{v_{1},v_{4}},f_{v_{2},v_{3}}).$$ Then $J_{2}(G_{1})={\rm rad}(I)$, since $f_{v_{1},v_{3}}^2 \in I$, $f_{v_{2},v_{6}}^2 \in I$, $f_{v_{1},v_{5}}^4 \in I$ and $f_{v_{3},v_{4}}^4 \in I$.\\ (2) $E(D)=\{\{v_{1},v_{2}\},\{v_{2},v_{3}\},\{v_{3},v_{4}\},\{v_{1},v_{4}\},\{v_{2},v_{4}\}\}$. Let $$I=(f_{v_{1},v_{5}}+f_{v_{2},v_{4}},f_{v_{2},v_{6}}+f_{v_{3},v_{4}},f_{v_{1},v_{2}},f_{v_{1},v_{4}},f_{v_{2},v_{3}}).$$ Then $J_{2}(G_{1})={\rm rad}(I)$, since $f_{v_{1},v_{5}}^2 \in I$, $f_{v_{2},v_{4}}^2 \in I$, $f_{v_{2},v_{6}}^4 \in I$ and $f_{v_{3},v_{4}}^4 \in I$.\\ Let  $G_{2}=T_{1}-\{v_1,v_{5}\}$ and $G_{3}=T_{2}-\{v_2,v_{6}\}$ be subgraphs of $G$. We have that ${\rm ara}(J_{2}(G_{2})) \leq r_{1}-1$ and ${\rm ara}(J_{2}(G_{3})) \leq r_{2}-1$. Thus ${\rm ara}(J_{2}(G)) \leq 5+\sum_{h=1}^{s}r_{h}-2$ and therefore ${\rm ara}(J_{2}(G))=3+\sum_{h=1}^{s}r_{h}$. \hfill $\square$

\begin{thm1} \label{12} Let $G$ be a connected graph on the vertex set $[n]$ with $n \geq 6$. Suppose that $G$ has a complete subgraph $\mathcal{K}_4$ and also $$G=\mathcal{K}_{4} \cup_{v_1} T_{1} \cup_{v_2} \cdots \cup_{v_s}T_{s},$$ where $\{v_{1},\ldots,v_{s}\} \subset V(\mathcal{K}_{4})$, $s \geq 2$, $v_{h}$ are pairwise distinct and $T_{h}$ are trees. Then ${\rm ara}(J_{2}(G))={\rm pd}(R/J_{2}(G))=n-1$. 
\end{thm1}

\noindent \textbf{Proof.} By \cite[Theorem 1.1]{EHH}, ${\rm depth}(R/J_{2}(G))=n+1$  and therefore ${\rm pd}(R/J_{2}(G))=n-1$. Thus ${\rm ara}(J_{2}(G)) \geq n-1$. Let $r_{h}$ be the number of edges of the tree $T_{h}$ for $1 \leq h \leq s$. The graph $G$ has $n=4+\sum_{h=1}^{s}r_{h}$ vertices. Let $\{v_{1},\ldots,v_{4}\}$ be the vertex set of $\mathcal{K}_4$ and let $G_{1}=\mathcal{K}_{4} \cup \{v_{1},v_{5}\} \cup \{v_{2},v_{6}\}$, $G_{2}=T_{1}-\{v_1,v_{5}\}$ and $G_{3}=T_{2}-\{v_2,v_{6}\}$ be subgraphs of $G$. Consider the ideal $$I=(f_{v_{1},v_{3}}+f_{v_{2},v_{6}},f_{v_{1},v_{5}}+f_{v_{3},v_{4}},f_{v_{1},v_{4}}+f_{v_{2},v_{3}},f_{v_{1},v_{2}},f_{v_{2},v_{4}}).$$ We have that $J_{2}(G_{1})={\rm rad}(I)$, since $f_{v_{1},v_{3}}^3 \in I$, $f_{v_{2},v_{6}}^3 \in I$, $f_{v_{1},v_{5}}^5 \in I$, $f_{v_{3},v_{4}}^5 \in I$, $f_{v_{1},v_{4}}^2 \in I$ and $f_{v_{2},v_{3}}^2 \in I$. Consequently ${\rm ara}(J_{2}(G_{1})) \leq 5$. Also ${\rm ara}(J_{2}(G_{2})) \leq r_{1}-1$ and ${\rm ara}(J_{2}(G_{3})) \leq r_{2}-1$. Thus ${\rm ara}(J_{2}(G)) \leq 5+\sum_{h=1}^{s}r_{h}-2$ and therefore ${\rm ara}(J_{2}(G))=3+\sum_{h=1}^{s}r_{h}$. \hfill $\square$

\begin{cor1} \label{2} Let $G$ be a connected graph on the vertex set $[n]$ with $n \geq 6$. Suppose that $G$ has a complete subgraph $\mathcal{K}_4$ and also $$G=\mathcal{K}_{4} \cup_{v_1} T_{1} \cup_{v_2} \cdots \cup_{v_s}T_{s},$$ where $\{v_{1},\ldots,v_{s}\} \subset V(\mathcal{K}_{4})$, $s \geq 2$, $v_{h}$ are pairwise distinct and $T_{h}$ are paths. Then $J_{2}(G)$ is a set-theoretic complete intersection. 
\end{cor1}

\noindent \textbf{Proof.} By \cite[Theorem 1.1]{EHH}, $J_{2}(G)$ is Cohen-Macaulay, so ${\rm dim}(R/J_{2}(G))=n+1$ and therefore ${\rm ht}(J_{2}(G))=n-1$. By Theorem \ref{12}, ${\rm ara}(J_{2}(G))={\rm ht}(J_{2}(G))$. \hfill $\square$

\begin{thm1} Let $G$ be a connected graph which is not a triangle and let $K$ be a field of any characteristic. Suppose that $J_{2}(G)$ is an almost complete intersection binomial edge ideal. Then $J_{2}(G)$ is a cohomologically complete intersection if and only if $J_{2}(G)$ is Cohen-Macaulay.
\end{thm1}
\noindent \textbf{Proof.} The implication $\Rightarrow$ is obvious from Proposition \ref{BasicProp}. Conversely suppose that $J_{2}(G)$ is Cohen-Macaulay. By \cite[Theorem 2]{Rin}, there are two cases. \begin{enumerate} \item G is a graph on the vertex set $\{u_{1},\ldots,u_{r},v_{1},\ldots,v_{s},w_{1},\ldots,w_{t}\}$ with $r \geq 2$, $s \geq 2$, $t \geq 2$ and edge set $$E(G)=\{\{u_{i},u_{i+1}\}: i=1,\ldots,r-1\} \cup \{\{v_{i},v_{i+1}\}:i=1,\ldots,s-1\} \cup$$ $ \ \ \ \ \ \ \ \ \ \cup \{\{w_{i},w_{i+1}\}:i=1,\ldots,t-1\} \cup \{\{u_{1},v_{1}\}, \{u_{1},w_{1}\},\{v_{1},w_{1}\}\}$.\\ Then ${\rm ht}(J_{2}(G))=r+s+t-1$. Let $S$ be the subgraph of $G$ with the edge set $E(S)=\{\{u_{1},v_{1}\}, \{u_{1},u_{2}\},\{v_{1},w_{1}\},\{u_{1},w_{1}\}\}$. Let $I=(f_{u_{1},u_{2}}+f_{v_{1},w_{1}},f_{u_{1},v_{1}},f_{u_{1},w_{1}})$, then $J_{2}(S)={\rm rad}(I)$ since $f_{u_{1},u_{2}}^2 \in I$ and $f_{v_{1},w_{1}}^2 \in I$. So ${\rm ara}(J_{2}(S)) \leq 3$ and also ${\rm ara}(J_{2}(G)) \leq 3+r+s+t-4$. Thus ${\rm ara}(J_{2}(G))=r+s+t-1$, so $J_{2}(G)$ is a set-theoretic complete intersection. \item G is a graph on the vertex set $\{u_{1},\ldots,u_{r},v_{1},\ldots,v_{s}\}$ with $r \geq 3$, $s \geq 3$ and edge set $$E(G)=\{\{u_{i},u_{i+1}\}: i=1,\ldots,r-1\} \cup \{\{v_{i},v_{i+1}\}:i=1,\ldots,s-1\} \cup$$ $\ \ \ \ \  \  \  \ \ \cup \{\{u_{1},v_{1}\}, \{u_{2},v_{2}\}\}$.\\ Then ${\rm ht}(J_{2}(G))=r+s-1$. Let $G_{1}$ be the subgraph of $G$ on the vertex set $\{u_{1},u_{2},u_{3},v_{1},v_{2},v_{3}\}$  with edges $\{u_{1},u_{2}\}$, $\{u_{2},v_{2}\}$, $\{v_{1},v_{2}\}$, $\{u_{1},v_{1}\}$, $\{u_{2},u_{3}\}$ and $\{v_{2},v_{3}\}$. Let $T_{1}=(\{u_{3},u_{4}\},\{u_{4},u_{5}\},\ldots,\{u_{r-1},u_{r}\})$ and $T_{2}=(\{v_{3},v_{4}\},\{v_{4},v_{5}\},\ldots,\{v_{s-1},v_{s}\})$. The ideals $J_{2}(T_{1})$ and $J_{2}(T_{2})$ are complete intersection, so ${\rm cd}(J_{2}(T_{1}))=r-3$ and ${\rm cd}(J_{2}(T_{2}))=s-3$. Then $G=G_{1} \cup_{u_3} T_{1} \cup_{v_3} T_{2}$. Let $G_2$ be the graph on the vertex set $\{u_{1},u_{2},u_{3},v_{1},v_{2},v_{3}\}$ with edges $\{u_{1},u_{2}\},$ $\{u_{2},v_{2}\}$, $\{v_{1},v_{2}\}$, $\{u_{1},v_{1}\}$, $\{u_{1},v_{2}\}$, $\{u_{2},u_{3}\}$ and $\{v_{2},v_{3}\}$. By Theorem \ref{11}, ${\rm cd}(J_{2}(G_{2}))=5$. Let $G_3$ be the graph on the vertex set $\{u_{1},u_{2},u_{3},v_{1},v_{2},v_{3}\}$  with edges $\{u_{1},u_{2}\}$, $\{u_{2},v_{2}\}$, $\{v_{1},v_{2}\}$, $\{u_{1},v_{1}\}$, $\{u_{2},v_{1}\}$, $\{u_{2},u_{3}\}$ and $\{v_{2},v_{3}\}$. By Theorem \ref{11}, ${\rm cd}(J_{2}(G_{3}))=5$.  Also $$J_{2}(G_{2})=P_{\emptyset}(G) \cap P_{\{u_2\}}(G) \cap P_{\{v_2\}}(G) \cap P_{\{u_{1},v_2\}}(G) \cap P_{\{u_{2},v_2\}}(G)$$ and $$J_{2}(G_{3})=P_{\emptyset}(G) \cap P_{\{u_2\}}(G) \cap P_{\{v_2\}}(G) \cap P_{\{u_{2},v_1\}}(G) \cap P_{\{u_{2},v_2\}}(G).$$ Thus $J_{2}(G_{1})=J_{2}(G_{2}) \cap J_{2}(G_{3})$, since $\mathcal{C}(G_{1})=\{\emptyset,\{u_{2}\},\{v_{2}\},\{u_{1},v_{2}\},\{u_{2},v_{1}\},$\\ $\{u_{2},v_{2}\}\}$. Let $G_{4}$ be the graph on the vertex set $\{u_{1},u_{2},u_{3},v_{1},v_{2},v_{3}\}$  with edges $\{u_{1},u_{2}\}$, $\{u_{2},v_{2}\}$, $\{v_{1},v_{2}\}$, $\{u_{1},v_{1}\}$, $\{u_{1},v_{2}\}$, $\{u_{2},v_{1}\}$, $\{u_{2},u_{3}\}$ and $\{v_{2},v_{3}\}$. By Theorem \ref{12}, ${\rm cd}(J_{2}(G_{4}))=5$. Notice that $J_{2}(G_{4})=J_{2}(G_2)+J_{2}(G_3)$. Since ${\rm cd}(J_{2}(G_{2}))=5$, ${\rm cd}(J_{2}(G_{3}))=5$ and ${\rm cd}(J_{2}(G_{4}))=5$, we deduce from the Mayer-Vietoris sequence  $$\cdots \rightarrow H_{J_{2}(G_{4})}^{i}(R) \rightarrow H_{J_{2}(G_{2})}^{i}(R) \oplus H_{J_{2}(G_{3})}^{i}(R) \rightarrow H_{J_{2}(G_{1})}^{i}(R) \rightarrow H_{J_{2}(G_{4})}^{i+1}(R) \rightarrow \cdots$$ that ${\rm cd}(J_{2}(G_{1})) \leq 5$. By \cite[Corollary 2.2]{DV}, ${\rm cd}(J_{2}(G)) \leq {\rm cd}(J_{2}(G_{1}))+{\rm cd}(J_{2}(T_1))+{\rm cd}(J_{2}(T_2))$, thus  ${\rm cd}(J_{2}(G)) \leq r+s-1$ and therefore ${\rm cd}(J_{2}(G))=r+s-1$. \hfill $\square$\\
\end{enumerate}

\begin{rem1} \label{almostbasic} {\rm If $\mu(J_{2}(G))={\rm ht}(J_{2}(G))+1$ and $J_{2}(G)$ is not Cohen-Macaulay, then ${\rm cd}(J_{2}(G))={\rm ht}(J_{2}(G))+1$, so ${\rm ara}(J_{2}(G))={\rm pd}(R/J_{2}(G))={\rm ht}(J_{2}(G))+1$.}
\end{rem1}

Let $G$ and $G'$ be two simple graphs with the vertex sets $[p]$ and $[q]$, respectively. The {\em join} of $G$ and $G'$, denoted by $G \ast G'$, is the graph with vertex set $[p] \cup [q]$ and the edge set $E(G \ast G')=E(G) \cup E(G') \cup \{\{i,j\}|i \in [p], j \in [q]\}$.

\begin{ex1} \label{Basicex} {\rm Let $G$ be the graph on the vertex set $[2]$ with only one edge $\{1,2\}$ and $G'$ be the graph on the vertex set $\{3,4,5\}$ with connected components $\{3,4\}$ and $\{5\}$. Then $G \ast G'$ is the graph on the vertex set $[5]$ with edges $ E(G \ast G')=\{\{1,2\},\{1,3\},\{1,4\},\{1,5\},\{2,3\},\{2,4\},\{2,5\},\{3,4\}\}$. From \cite[Theorem 3.9]{KuSa} we have that ${\rm depth}(R/J_{2}(G \ast G'))=5$ and hence by Auslander-Buchsbaum formula, ${\rm pd}(R/J_{2}(G \ast G'))=5$. Thus ${\rm ara}(J_{2}(G \ast G')) \geq 5$. Let $$I=(f_{1,3}+f_{2,5},f_{1,4}+f_{2,3},f_{1,5}+f_{3,4},f_{1,2},f_{2,4}).$$ Then
$J_{2}(G \ast G')={\rm rad}(I)$, since $f_{1,3}^3 \in I$, $f_{2,5}^3 \in I$, $f_{1,4}^2 \in I$, $f_{2,3}^2 \in I$,  $f_{1,5}^5 \in I$ and $f_{3,4}^5 \in I$. So ${\rm ara}(J_{2}(G \ast G')) \leq 5$ and therefore ${\rm ara}(J_{2}(G \ast G'))=5$.}
\end{ex1}

\begin{prop1} \label{join} Let $G$ be a graph on the vertex set $[p]$ with $r \geq 1$ edges and $G'$ be a graph on the vertex set $[q]$ with $t \geq 1$ edges, where $q \geq 3$. Then $${\rm ara}(J_{2}(G \ast G'))  \leq pq+r+t-3.$$
\end{prop1}

\noindent \textbf{Proof.} Suppose that $\{1,2\}$ and $\{3,4\}$ are edges of $G$ and $G'$, respectively. Also assume that $5$ is a vertex of $G'$. Consider the subgraph $S$ of $G \ast G'$ with the edge set $ E(S)=\{\{1,2\},\{1,3\},\{1,4\},\{1,5\},\{2,3\},\{2,4\},\{2,5\},\{3,4\}\}$. By Example \ref{Basicex}, ${\rm ara}(J_{2}(S))=5$. Since the graph $G \ast G'$ has $pq+r+t$ edges, we get ${\rm ara}(J_{2}(G \ast G'))  \leq pq+r+t-3$. \hfill $\square$\\

We concentrate now on graphs of the form $G \ast 2 \mathcal{K}_{1}$, where $G$ is a graph on the vertex set $[p]$ and $2 \mathcal{K}_{1}$ is the graph consisting of two isolated vertices. By \cite[Theorem 5.3]{RSK}, ${\rm depth}(R/J_{2}(G \ast 2 \mathcal{K}_{1}))=4$, so using the Auslander-Buchsbaum formula we get ${\rm pd}(R/J_{2}(G \ast 2 \mathcal{K}_{1}))=2p$. Thus $2p \leq {\rm ara}(J_{2}(G \ast 2 \mathcal{K}_{1}))$. 

A {\em null graph} is a graph in which there are no edges between its vertices.

\begin{prop1} \label{join} Let $G$ be a null graph on the vertex set $[p]$ where $p \geq 2$. Then ${\rm ara}(J_{2}(G \ast 2 \mathcal{K}_{1}))=2p$.
\end{prop1}
\noindent \textbf{Proof.} The graph $G \ast 2 \mathcal{K}_{1}$ is complete bipartite with $p+2$ vertices. Also $\mu(J_{2}(G \ast 2 \mathcal{K}_{1}))=2p$, so ${\rm ara}(J_{2}(G \ast 2 \mathcal{K}_{1})) \leq 2p$ and therefore ${\rm ara}(J_{2}(G \ast 2 \mathcal{K}_{1}))=2p$. \hfill $\square$

\begin{rem1} \label{completebipartite} {\rm For the complete bipartite graph $\mathcal{K}_{2,p}$ with $p+2$ vertices we have that ${\rm ara}(J_{2}(\mathcal{K}_{2,p}))=2p$.}
\end{rem1}

\begin{thm1} \label{even} For every even positive integer $n$ there is a graph $G$ such that ${\rm depth}(R/J_{2}(G))=4$ and ${\rm ara}(J_{2}(G))={\rm pd}(R/J_{2}(G))=n$.
\end{thm1}

\noindent \textbf{Proof.} Let $n=2p$, where $p \geq 1$ is an integer. Consider the complete bipartite graph $\mathcal{K}_{2,p}$ with $p+2$ vertices. By Remark \ref{completebipartite},  ${\rm ara}(J_{2}(\mathcal{K}_{2,p}))=2p$. \hfill $\square$\\

Given a connected and simple graph $G$ on the vertex set $[n]$ with edge set $E(G)=\{e_{1},\ldots,e_{r}\}$, we write $\{e_{1},\ldots,e_{r}\} \ast 2 \mathcal{K}_{1}:=G \ast 2 \mathcal{K}_{1}$.

\begin{ex1} \label{Examplejoin} {\rm Let $G$ be a graph on the vertex set $[p]$ with exactly one edge $e=\{1,2\}$, where $p \geq 2$. Then $J_{2}(G \ast 2 \mathcal{K}_{1})=J_{2}(\{e\} \ast 2 \mathcal{K}_{1})+J_{2}(\mathcal{K}_{2,p-2})$, where $\mathcal{K}_{2,p-2}$ is a complete bipartite graph. Notice that $J_{2}(\mathcal{K}_{2,p-2})=(0)$ if $p=2$. Let $\{i,j\}$ be the vertex set of $2 \mathcal{K}_{1}$ and consider the ideal $I=(f_{1,i}+f_{2,j},f_{1,2},f_{1,j},f_{2,i})$. Then $J_{2}(\{e\}  \ast 2 \mathcal{K}_{1})={\rm rad}(I)$ since $f_{1,i}^2 \in I$ and $f_{2,j}^2 \in I$. Thus ${\rm ara}(J_{2}(\{e\}  \ast 2 \mathcal{K}_{1})) \leq 4$, so ${\rm ara}(J_{2}(G \ast 2 \mathcal{K}_{1})) \leq 2p$ since ${\rm ara}(J_{2}(\mathcal{K}_{2,p-2}))=2p-4$. Consequently ${\rm ara}(J_{2}(G \ast 2 \mathcal{K}_{1})) =2p$.}
\end{ex1}

\begin{thm1} \label{An} Let $G$ be a graph on the vertex set $[p]$ which has at most $3$ edges, where $p \geq 2$. Then ${\rm ara}(J_{2}(G \ast 2 \mathcal{K}_{1}))={\rm pd}(R/J_{2}(G \ast 2 \mathcal{K}_{1}))=2p$.
\end{thm1}
\noindent \textbf{Proof.} It is enough to show that  ${\rm ara}(J_{2}(G \ast 2 \mathcal{K}_{1})) \leq 2p$. Let $\{i,j\}$ be the vertex set of $2\mathcal{K}_{1}$. We distinguish the following cases:
\begin{enumerate} \item $G$ is a null graph. By Proposition \ref{join}, ${\rm ara}(J_{2}(G \ast 2 \mathcal{K}_{1}))=2p$. 
\item $G$ has exactly one edge $e$. By Example \ref{Examplejoin}, ${\rm ara}(J_{2}(G \ast 2 \mathcal{K}_{1})) =2p$.
\item $G$ has exactly two edges $e_{1}$ and $e_{2}$. There are two cases: (i) $e_{1} \cap e_{2}=\emptyset$. Then $J_{2}(G \ast 2 \mathcal{K}_{1})=J_{2}(\{e_1\} \ast 2 \mathcal{K}_{1})+J_{2}(\{e_2\} \ast 2 \mathcal{K}_{1})+J_{2}(\mathcal{K}_{2,p-4})$. Also ${\rm ara}(J_{2}(\{e_1\} \ast 2 \mathcal{K}_{1}))=4$, ${\rm ara}(J_{2}(\{e_2\} \ast 2 \mathcal{K}_{1}))=4$ and ${\rm ara}(J_{2}(\mathcal{K}_{2,p-4}))=2p-8$. Thus ${\rm ara}(J_{2}(G \ast 2 \mathcal{K}_{1})) \leq 2p$.\\ (ii) $e_{1} \cap e_{2} \neq \emptyset$. Let $e_{1}=\{1,2\}$ and $e_{2}=\{2,3\}$. Then $J_{2}(G \ast 2 \mathcal{K}_{1})=J_{2}(\{e_1,e_{2}\} \ast 2 \mathcal{K}_{1})+J_{2}(\mathcal{K}_{2,p-3})$. Let $I=(f_{1,2},f_{2,3},f_{1,j},f_{3,i},f_{1,i}+f_{2,j},f_{2,i}+f_{3,j})$. We have that $J_{2}(\{e_1,e_{2}\} \ast 2 \mathcal{K}_{1})={\rm rad}(I)$, since $f_{1,i}^2 \in I$, $f_{2,j}^2 \in I$, $f_{2,i}^2 \in I$ and $f_{3,j}^2 \in I$. So ${\rm ara}(J_{2}(\{e_1,e_{2}\} \ast 2 \mathcal{K}_{1})) \leq 6$ and therefore ${\rm ara}(J_{2}(G \ast 2 \mathcal{K}_{1})) \leq 2p$, since ${\rm ara}(J_{2}(\mathcal{K}_{2,p-3}))=2p-6$.
\item $G$ has exactly three edges $e_{1}$, $e_{2}$ and $e_{3}$. There are five cases: (i) $e_{r} \cap e_{s} = \emptyset$, for every $1 \leq r<s \leq 3$. Then $J_{2}(G \ast 2 \mathcal{K}_{1})=J_{2}(\{e_1\} \ast 2 \mathcal{K}_{1})+J_{2}(\{e_2\} \ast 2 \mathcal{K}_{1})+J_{2}(\{e_3\} \ast 2 \mathcal{K}_{1})+J_{2}(\mathcal{K}_{2,p-6})$. Also ${\rm ara}(J_{2}(\{e_i\} \ast 2 \mathcal{K}_{1}))=4$, for $1 \leq i \leq 3$, and ${\rm ara}(J_{2}(\mathcal{K}_{2,p-6}))=2p-12$. Thus ${\rm ara}(J_{2}(G \ast 2 \mathcal{K}_{1})) \leq 2p$.\\
(ii) $e_{1} \cap e_{2} \neq \emptyset$, $e_{1} \cap e_{3}=\emptyset$ and $e_{2} \cap e_{3}=\emptyset$. Then $J_{2}(G \ast 2 \mathcal{K}_{1})=J_{2}(\{e_1,e_{2}\} \ast 2 \mathcal{K}_{1})+J_{2}(\{e_3\} \ast 2 \mathcal{K}_{1})+J_{2}(\mathcal{K}_{2,p-5})$.  Also ${\rm ara}(J_{2}(\{e_1,e_{2}\} \ast 2 \mathcal{K}_{1}))=6$, ${\rm ara}(J_{2}(\{e_3\} \ast 2 \mathcal{K}_{1}))=4$ and ${\rm ara}(J_{2}(\mathcal{K}_{2,p-5}))=2p-10$. So ${\rm ara}(J_{2}(G \ast 2 \mathcal{K}_{1})) \leq 2p$.\\
(iii) $\{e_{1},e_{2},e_{3}\}$ is a path in $G$. Let $e_{1}=\{1,2\}$, $e_{2}=\{2,3\}$ and $e_{3}=\{3,4\}$. Then $J_{2}(G \ast 2 \mathcal{K}_{1})=J_{2}(\{e_{1},e_{2},e_{3}\} \ast 2 \mathcal{K}_{1})+J_{2}(\mathcal{K}_{2,p-4})$. Let $I=(f_{1,2},f_{2,3},f_{3,4},f_{1,j},f_{4,i},f_{1,i}+f_{2,j},f_{2,i}+f_{3,j},f_{3,i}+f_{4,j})$. We have that $J_{2}(\{e_{1},e_{2},e_{3}\} \ast 2 \mathcal{K}_{1})={\rm rad}(I)$, since $f_{1,i}^2 \in I$, $f_{2,j}^2 \in I$, $f_{2,i}^4 \in I$, $f_{3,j}^4 \in I$, $f_{3,i}^2 \in I$ and $f_{4,j}^2 \in I$. So ${\rm ara}(J_{2}(\{e_1,e_{2},e_{3}\} \ast 2 \mathcal{K}_{1})) \leq 8$ and therefore ${\rm ara}(J_{2}(G \ast 2 \mathcal{K}_{1})) \leq 2p$, since ${\rm ara}(J_{2}(\mathcal{K}_{2,p-4}))=2p-8$.\\
(iv) $e_{1} \cap e_{2} \cap e_{3} \neq \emptyset$. Let $e_{1}=\{1,2\}$, $e_{2}=\{1,3\}$ and $e_{3}=\{1,4\}$. Then $J_{2}(G \ast 2 \mathcal{K}_{1})=J_{2}(\{e_{1},e_{2},e_{3}\} \ast 2 \mathcal{K}_{1})+J_{2}(\mathcal{K}_{2,p-4})$. Let $I=(f_{1,2}+f_{3,j},f_{1,i}+f_{2,j},f_{1,j}+f_{4,i},f_{1,3},f_{1,4},f_{2,i},f_{3,i},f_{4,j})$. We have that $J_{2}(\{e_{1},e_{2},e_{3}\} \ast 2 \mathcal{K}_{1})={\rm rad}(I)$, since $f_{1,2}^2 \in I$, $f_{3,j}^2 \in I$, $f_{1,i}^3 \in I$, $f_{2,j}^3 \in I$, $f_{1,j}^2 \in I$ and $f_{4,i}^2 \in I$. So ${\rm ara}(J_{2}(\{e_1,e_{2},e_{3}\} \ast 2 \mathcal{K}_{1})) \leq 8$ and therefore ${\rm ara}(J_{2}(G \ast 2 \mathcal{K}_{1})) \leq 2p$, since ${\rm ara}(J_{2}(\mathcal{K}_{2,p-4}))=2p-8$.\\
(v) $\{e_{1},e_{2},e_{3}\}$ is the edge set of a triangle. Let $e_{1}=\{1,2\}$, $e_{2}=\{2,3\}$ and $e_{3}=\{1,3\}$. Then $J_{2}(G \ast 2 \mathcal{K}_{1})=J_{2}(\{e_{1},e_{2},e_{3}\} \ast 2 \mathcal{K}_{1})+J_{2}(\mathcal{K}_{2,p-3})$. Let $I=(f_{1,2}+f_{3,i},f_{1,3}+f_{2,j},f_{1,j}+f_{2,i},f_{1,i},f_{2,3},f_{3,j})$. We have that $J_{2}(\{e_{1},e_{2},e_{3}\} \ast 2 \mathcal{K}_{1})={\rm rad}(I)$, since $f_{1,2}^3 \in I$, $f_{3,i}^3 \in I$, $f_{1,3}^2 \in I$, $f_{2,j}^2 \in I$, $f_{1,j}^5 \in I$ and $f_{2,i}^5 \in I$. So ${\rm ara}(J_{2}(\{e_1,e_{2},e_{3}\} \ast 2 \mathcal{K}_{1})) \leq 6$ and therefore ${\rm ara}(J_{2}(G \ast 2 \mathcal{K}_{1})) \leq 2p$, since ${\rm ara}(J_{2}(\mathcal{K}_{2,p-3}))=2p-6$. \hfill $\square$

\end{enumerate}

\end{document}